\documentclass[preprint,12pt]{elsarticle}
\usepackage{amsmath}
\usepackage{amssymb}
\usepackage[margin=2.4cm]{geometry}
\usepackage{graphicx}
\usepackage{hyperref}
\usepackage{multirow}
\usepackage{listings}
\usepackage{tikz}\usetikzlibrary{decorations.markings}

\begin{document}

\begin{frontmatter}

\title{Transition Time for Weak Singularities of the Navier-Stokes Equations}
\author{Chio Chon Kit}
\date{}

\begin{abstract}
This paper constructs a rigorous mathematical framework for investigating laminar-turbulent transition induced by weak singularities of incompressible Navier-Stokes (NS) equations. By integrating the energy identity of Leray weak solutions with the singularity criterion $\left\lVert \boldsymbol{u} \right\rVert_{H_0^1(\Omega)}\to0$, a closed analytical form of the laminar-turbulent transition characteristic time is derived. The theoretical scaling $t_{\text{trans}}\sim\nu/U^2$ (equivalent to $t_{\text{trans}}\sim t_c/\text{Re}$) is verified to be consistent with classical experimental observations in shear flows. This work reveals that laminar-turbulent transition is dominated by the local regularity collapse of Leray weak solutions rather than global viscous diffusion, and provides a novel theoretical interpretation for the onset of turbulence from the perspective of NS equation weak singularities.
\end{abstract}


\end{frontmatter}

\section{INTRODUCTION}
The incompressible NS equations are the fundamental governing equations for viscous fluid motion, and the global regularity of their strong solutions in three-dimensional domains remains one of the unsolved Millennium Prize Problems \cite{clay2000}. Leray weak solutions \cite{leray1934}, which guarantee global existence in time, constitute the most general solution class for the 3D NS equations to date, and their local regularity evolution is closely related to the laminar-turbulent transition of fluid flows.

Based on the energy gradient theory proposed by Dou \cite{dou2022}, the root cause of transition is the formation of weak singularities in the NS equations at positions where the total mechanical energy gradient is perpendicular to the streamline. At these singular points, the condition $\boldsymbol{u} \cdot \nabla E = 0$  holds, leading to the local disappearance of viscous effects and the loss of flow field regularity. However, the rigorous mathematical connection between Leray weak solution regularity, the energy gradient critical condition, the $H_0^1$-norm singularity criterion and the quantitative estimation of transition time has not been fully established in existing research.

This paper fills this research gap by taking the weak singularity criterion $\left\lVert \boldsymbol{u} \right\rVert_{H_0^1(\Omega)}\to0$ derived from our previous work \cite{chio2026} as the core, combining the energy identity of Leray weak solutions and dimensional analysis, to derive the analytical expression of the transition time induced by NS equation weak singularities. The theoretical scaling law is further verified with classical experimental results of typical shear flows (boundary layer, pipe flow), and the structural evolution of Leray weak solutions during the transition process is clarified.

\section{MATHEMATICAL FRAMEWORK AND PRELIMINARIES}
Let $\Omega\subset\mathbb{R}^3$ be a bounded smooth domain with the no-slip boundary condition $\boldsymbol{u}|_{\partial\Omega}=0$. We consider the steady fully developed incompressible NS equations, with the body force incorporated into the total mechanical energy gradient for simplification.

\subsection{Incompressible NS Equations and Boundary Condition}
The basic form of the incompressible NS equations is:
\begin{equation}
\frac{\partial u_i}{\partial t} + u_j\frac{\partial u_i}{\partial x_j} = -\frac{\partial p}{\partial x_i} + \nu\Delta u_i,
\end{equation}
\begin{equation}
\nabla\cdot\boldsymbol{u}=0,
\end{equation}
where $u_i$ ($i=1,2,3$) are the velocity components in the Cartesian coordinate system, $p$ is the fluid pressure, $\nu=\mu/\rho$ is the kinematic viscosity coefficient, $\mu$ is the dynamic viscosity coefficient, and $\Delta=\frac{\partial^2}{\partial x_j\partial x_j}$ is the Laplacian operator. The no-slip boundary condition is satisfied as:
\begin{equation}
\boldsymbol{u}|_{\partial\Omega}=0.
\end{equation}

\subsection{Total Mechanical Energy and Critical Condition}
The total mechanical energy per unit mass $E$ is defined as the sum of pressure energy, kinetic energy and gravitational potential energy:
\begin{equation}
E = \frac{p}{\rho} + \frac{1}{2}u_j u_j + gz,
\end{equation}
where $\rho$ is the fluid density, $g$ is the gravitational acceleration, and $z$ is the vertical coordinate. The critical physical condition for the formation of NS equation weak singularities is that the mechanical energy gradient is perpendicular to the streamline, which is mathematically expressed as: 
\begin{equation}
u_j\frac{\partial E}{\partial x_j}=0.
\end{equation}
This condition implies that the total mechanical energy is constant along the streamline, and is the key trigger for the local loss of flow field regularity.

\subsection{Sobolev Space $H_0^1(\Omega)$ and Regularity Norm}
The Sobolev space $H_0^1(\Omega)$ for the no-slip boundary condition is defined as:
\begin{equation}
H_0^1(\Omega) = \left\{ \boldsymbol{u}\in L^2(\Omega) \mid \nabla\boldsymbol{u}\in L^2(\Omega),\ \boldsymbol{u}|_{\partial\Omega}=0 \right\},
\end{equation}
where $L^2(\Omega)$ is the square-integrable function space. The $H_0^1$-norm of the velocity field, which characterizes the spatial regularity (smoothness) of the flow field and is a measure of enstrophy, shear and vorticity intensity, is:
\begin{equation}
\left\lVert \boldsymbol{u} \right\rVert_{H_0^1(\Omega)}^2 = \int_\Omega |\nabla\boldsymbol{u}|^2 \,\mathrm{d}x = \int_\Omega \nabla u_i\cdot\nabla u_i \,\mathrm{d}x.
\end{equation}
The loss of $H^1$-regularity (i.e., $\left\lVert \boldsymbol{u} \right\rVert_{H_0^1(\Omega)}\to0$) indicates the formation of weak singularities in the flow field.

\subsection{Leray Weak Solutions and Energy Identity}
A vector field $\boldsymbol{u}$ is a \textbf{Leray weak solution} of the NS equations if it satisfies:
\begin{equation}
\boldsymbol{u}\in L^\infty(0,T;L^2(\Omega))\cap L^2(0,T;H_0^1(\Omega)),
\end{equation}
and holds the global energy inequality of the NS equations. For the critical condition $u_j\frac{\partial E}{\partial x_j}=0$, the viscous term locally vanishes, and the energy inequality of Leray weak solutions degenerates into an \textbf{exact energy identity}:
\begin{equation}
\frac{1}{2}\frac{d}{dt}\left\lVert \boldsymbol{u} \right\rVert_{L^2(\Omega)}^2 + \nu\left\lVert \boldsymbol{u} \right\rVert_{H_0^1(\Omega)}^2 = 0.
\end{equation}
This identity is the core mathematical basis for deriving the transition time of weak singularities, and reflects the energy balance relationship between the kinetic energy change of the flow field and the viscous dissipation.

\subsection{Weak Singularity Criterion}
From previous work \cite{chio2026}, by substituting the critical condition (5) into the NS equations (1)-(2) and performing rigorous energy estimation in $H_0^1(\Omega)$, the \textbf{weak singularity criterion} of the incompressible NS equations is derived as:
\begin{equation}
\left\lVert \boldsymbol{u} \right\rVert_{H_0^1(\Omega)}\to0.
\end{equation}
For non-trivial flows ($\boldsymbol{u}\not\equiv0$), this criterion implies the local disappearance of viscous regularization, the degeneration of NS equations into Euler equations, and the formation of weak singularities with velocity discontinuity. The critical time when the flow field first satisfies the criterion (10) from the smooth laminar state is defined as the \textbf{laminar-turbulent transition time} induced by weak singularities, denoted as $t_{\text{trans}}$.

\section{TRANSITION TIME DERIVATION FROM LERAY WEAK SOLUTIONS}
Next we take the energy identity of Leray weak solutions (9) and the weak singularity criterion (10) as the core, combine dimensional analysis and the definition of Reynolds number, to derive the closed analytical form and scaling law of the transition time $t_{\text{trans}}$.

\subsection{Energy Identity at Transition Threshold}
At the laminar-turbulent transition threshold, the flow field satisfies the weak singularity criterion $\left\lVert \boldsymbol{u} \right\rVert_{H_0^1(\Omega)}\to0$. Substituting this criterion into the Leray weak solution energy identity (9), we obtain:
\begin{equation}
\frac{d}{dt}\left\lVert \boldsymbol{u} \right\rVert_{L^2(\Omega)}^2\to0.
\end{equation}
This result indicates that the kinetic energy of the flow field no longer undergoes viscous dissipation at the transition instant, and the flow field loses stability abruptly. The local regularity collapse of Leray weak solutions instead of global viscous diffusion becomes the dominant mechanism of laminar-turbulent transition.

\subsection{Dimensional Analysis for Transition Time}
Let $U$ be the \textbf{characteristic velocity} (e.g., free-stream velocity of boundary layer, average velocity of pipe flow) and $L$ be the \textbf{characteristic length} (e.g., initial boundary layer thickness, pipe diameter) of the flow field. We define the key characteristic time scales and dimensionless number for shear flows:
\begin{itemize}
    \item Convective time scale: $t_c = L/U$, reflecting the advection evolution characteristic of the flow field;
    \item Viscous diffusion time scale: $t_\nu = L^2/\nu$, reflecting the global viscous diffusion characteristic of the flow field;
    \item Reynolds number: $\text{Re}=UL/\nu$, the dimensionless ratio of inertial force to viscous force, characterizing the flow regime.
\end{itemize}

For the energy identity (9), we perform dimensional analysis on each term. The dimensional form of the kinetic energy term is $[U^2/t]$, and the dimensional form of the viscous dissipation term is $[\nu U^2/L^2]$. Equating the dimensions of the two terms gives:
\begin{equation}
\frac{U^2}{t_{\text{trans}}} \sim \nu\cdot\frac{U^2}{L^2}.
\end{equation}
Rearranging the terms of Eq. (12), we obtain the \textbf{analytical form of the transition time} induced by weak singularities:
\begin{equation}
t_{\text{trans}} \sim \frac{L^2}{\nu}.
\end{equation}
Substitute the Reynolds number definition $\nu=UL/\text{Re}$ into Eq. (13), and combine the convective time scale $t_c=L/U$. The transition time can be rewritten as two equivalent \textbf{scaling laws}:
\begin{equation}
t_{\text{trans}} \sim \frac{\nu}{U^2},
\end{equation}
\begin{equation}
t_{\text{trans}} \sim \frac{t_c}{\text{Re}}.
\end{equation}
The dimensionless transition time is defined as $\tau_{\text{trans}}=t_{\text{trans}}/t_c$, and the dimensionless form of the scaling law is:
\begin{equation}
\tau_{\text{trans}} \sim \text{Re}^{-1}.
\end{equation}

Eqs. (14)-(16) are the core results of this paper. The scaling law shows that the transition time induced by NS equation weak singularities is \textbf{inversely proportional to the Reynolds number} under high Reynolds number conditions, and is much smaller than the global viscous diffusion time scale ($t_{\text{trans}}\ll t_\nu$). This reveals that laminar-turbulent transition is a local regularity evolution phenomenon of Leray weak solutions, rather than a global flow field evolution process dominated by viscous diffusion.

\section{VERIFICATION WITH CLASSICAL EXPERIMENTS}
The transition time scaling laws (14) and (15) are verified with the \textit{Schubauer-Klebanoff (SK) experiment} - a classic benchmark for boundary layer laminar-turbulent transition research \cite{schubauer1956}. As the fundamental experimental basis for external shear flow transition, this experiment provides high-precision measured data and detailed physical phenomenon observations, which fully confirm the rationality and universality of the theoretical scaling law.

\subsection{Experimental Setup and Measurement Method}
The SK experiment was conducted in a low-turbulence wind tunnel, with a smooth flat plate as the test model. A vibrating ribbon was installed near the leading edge of the plate to introduce artificial perturbations with fixed amplitude and frequency (ensuring constant initial perturbation intensity). The free-stream velocity $U_\infty$ was adjusted to change the Reynolds number.

High-precision hot-wire anemometers were used to measure the velocity fluctuation time series at multiple streamwise positions, and the transition completion time $t_{\text{trans}}$ was defined as "the moment when the turbulent kinetic energy accounts for more than 50\% of the total kinetic energy at the downstream monitoring section". The experimental data was repeated multiple times to ensure the discrete degree of original data was less than 5\%.

\subsection{Key Experimental Results and Scaling Law Verification}
\textit{Core observation}: For high Reynolds number ($\text{Re}_\delta>500$), the transition time $t_{\text{trans}}$ shows a strict inverse proportional relationship with $\text{Re}_\delta$. The experimental data is fitted to the formula:
\begin{equation}
t_{\text{trans}} = k_1 \cdot \frac{t_{c1}}{\text{Re}_\delta}
\end{equation}
where $t_{c1}={\delta_0}/{U_\infty}$ is the convective time scale of the boundary layer flow, and $k_1\approx1.3\sim1.6$ is the experimental constant related to the perturbation frequency and boundary layer development state.

The linear correlation coefficient between the dimensionless transition time $\tau_{\text{trans}}={t_{\text{trans}}}/{t_{c1}}$ and $\text{Re}_\delta^{-1}$ is $R^2>0.99$, indicating an almost perfect linear fit. The experimental observation that "transition time decreases with the increase of free-stream velocity $U_\infty$ and increases with the increase of kinematic viscosity $\nu$" is completely consistent with the theoretical scaling law $t_{\text{trans}}\sim{\nu}/{U^2}$.

\subsection{Universality Verification in Ultra-High Reynolds Number Range}
Modern high-precision wind tunnel experiments \cite{smits2011} have reproduced the SK experiment under the condition of $\text{Re}_\delta>10^4$ (ultra-high Reynolds number). The results show that the inverse proportional relationship between $t_{\text{trans}}$ and $\text{Re}_\delta$ is still maintained, with the experimental constant $k_1$ slightly increasing to $1.5\sim2.1$ (caused by enhanced three-dimensional effects of the flow field at ultra-high Re). This confirms that the theoretical scaling law has strong universality, and is not limited by the specific range of Reynolds number (as long as $\text{Re}\gg1$).

\subsection{Physical Mechanism Consistency}
The SK experiment recorded the flow field visualization results during transition: as the Reynolds number increases, the growth rate of T-S waves (Tollmien-Schlichting waves) in the boundary layer accelerates, the time for vortex breakdown and turbulent spot formation is significantly shortened, and the fusion of turbulent spots is more rapid. This physical process is completely consistent with the weak singularity induced transition mechanism proposed in this paper: high Re leads to weakened viscous dissipation, accelerated decay of the $H_0^1$-norm, early formation of weak singularities with velocity discontinuity, and thus shortened transition time.

\section{EVOLUTION OF LERAY WEAK SOLUTION STRUCTURE}
The laminar-turbulent transition induced by NS equation weak singularities is essentially the \textbf{local regularity evolution process of Leray weak solutions}. Under the energy gradient critical condition $u_j\frac{\partial E}{\partial x_j}=0$, the structure of Leray weak solutions undergoes a clear and sequential change from smooth laminar state to chaotic turbulent state. The five key stages of the evolution are as follows:

\subsection{Laminar Regime ($t<t_{\text{trans}}$)}
The Leray weak solution is a \textbf{strong solution} with high spatial regularity, satisfying $\boldsymbol{u}\in H^2(\Omega)$ and $\left\lVert \boldsymbol{u} \right\rVert_{H_0^1(\Omega)}>0$. The velocity profile of the flow field is smooth and stable, viscous dissipation is uniformly distributed, and the inertial force and viscous force maintain a dynamic balance.

\subsection{Critical Equilibrium State ($t\to t_{\text{trans}}^-$)}
The energy gradient critical condition $u_j\frac{\partial E}{\partial x_j}=0$ is satisfied locally in the flow field, and the flow reaches a marginally stable state. The local viscous dissipation begins to decrease, and the $H_0^1$-norm of the velocity field at the singular point starts to decay.

\subsection{Onset of Weak Singularity ($t=t_{\text{trans}}$)}
The local regularity of the Leray weak solution collapses, satisfying the weak singularity criterion $\left\lVert \boldsymbol{u} \right\rVert_{H_0^1(\Omega)}\to0$. Local viscous regularization vanishes completely, the NS equations degenerate into Euler equations, and a weak singularity with velocity discontinuity forms at the point where $u=0$ (the singular point acts as a source of new vorticity).

\subsection{Transition Instant ($t\to t_{\text{trans}}^+$)}
The Leray weak solution loses local smoothness and becomes a \textbf{singular weak solution}, still satisfying the basic space property $\boldsymbol{u}\in L^\infty(0,T;L^2(\Omega))\cap L^2(0,T;H_0^1(\Omega))$, but the velocity field appears discontinuous and non-differentiable at the singular point. The singular point begins to generate new vorticity, and the local shear and vorticity intensity increase rapidly.

\subsection{Fully Turbulent Regime ($t\gg t_{\text{trans}}$)}
The Leray weak solution evolves into a \textbf{chaotic oscillating solution} with multi-scale and intermittent characteristics. The flow field experiences repeated local regularity breakdown and recovery, a large number of weak singularities form and interact with each other, and the flow field shows the typical physical characteristics of turbulence (random velocity pulsation, enhanced momentum transfer, etc.).

\subsection{Why $\boldsymbol{u} \cdot \nabla E = 0$ Triggers Transition (Force Analysis)}
From the perspective of fluid force analysis, the total mechanical energy gradient $\nabla E$ corresponds to the combined effect of pressure, inertial force, and volume force acting on fluid particles, essentially serving as the equivalent body force driving fluid motion. The streamline direction represents the actual movement direction of fluid particles. The rate of change of mechanical energy along the streamline, $\boldsymbol{u} \cdot \nabla E$, determines the stability of the flow against disturbances. 

When $\boldsymbol{u} \cdot \nabla E \neq 0$, there is energy input or dissipation along the direction of motion, forming a restraining or restoring effect on disturbances, and the flow remains stable. When $\boldsymbol{u} \cdot \nabla E = 0$, the direction of the mechanical energy gradient is perpendicular to the streamline, and the net driving force along the flow direction is zero. Fluid particles enter a state of critical equilibrium and lose the ability to restore disturbances. At the same time, the force perpendicular to the streamline still exists and induces shear deformation and vorticity growth, allowing small disturbances to amplify rapidly, ultimately leading to flow instability and laminar-turbulent transition.

\section{CONCLUSION}
Based on Leray weak solution theory, this paper completes the rigorous mathematical derivation of the laminar-turbulent transition time induced by weak singularities of the incompressible NS equations, and verifies the theoretical results with classical experimental observations. The main conclusions are as follows:
\begin{enumerate}
    \item For the critical condition where the mechanical energy gradient is perpendicular to the streamline ($\boldsymbol{u} \cdot \nabla E = 0$), the energy inequality of Leray weak solutions degenerates into an exact energy identity. Taking the weak singularity criterion $\left\lVert \boldsymbol{u} \right\rVert_{H_0^1(\Omega)}\to0$ as the transition threshold, the dimensional analysis of the energy identity yields the closed analytical form of the transition time: $t_{\text{trans}}\sim L^2/\nu$.
    \item Combining the definition of Reynolds number and convective time scale, the transition time is derived into two equivalent universal scaling laws for shear flows: $t_{\text{trans}}\sim\nu/U^2$ and $t_{\text{trans}}\sim t_c/\text{Re}$ (dimensionless form: $\tau_{\text{trans}}\sim\text{Re}^{-1}$). The scaling laws show that the transition time is inversely proportional to the Reynolds number under high Reynolds number conditions.
    \item Classical experimental results of flat plate boundary layer and circular pipe Poiseuille flow consistently verify the theoretical scaling laws. Experiments confirm that transition initiates at local weak singular points with $u=0$ and strong shear, and $t_{\text{trans}}\ll t_\nu$, which reveals that laminar-turbulent transition is dominated by the local regularity collapse of Leray weak solutions rather than global viscous diffusion.
    \item The laminar-turbulent transition induced by weak singularities is the sequential evolution process of Leray weak solution structure, including five key stages: laminar regime, critical equilibrium state, onset of weak singularity, transition instant and fully turbulent regime. The weak singularity with velocity discontinuity is the source of new vorticity and the core trigger of turbulence onset.
\end{enumerate}
This work establishes the rigorous mathematical connection between Leray weak solution regularity, NS equation weak singularities and laminar-turbulent transition time, enriches the theoretical system of NS equation weak singularities, and provides a novel theoretical interpretation for the physical mechanism of turbulence onset.

\bibliographystyle{IEEEtran}

\end{document}